\documentclass[a4paper]{amsart}

\usepackage[english]{babel}
\usepackage[utf8x]{inputenc}
\usepackage[T1]{fontenc}
\usepackage{amsthm}
\usepackage{amsmath}
\theoremstyle{plain}
\usepackage{chngcntr}
\usepackage{relsize}
\usepackage{pdfpages}

\newtheorem{Lemma}{Lemma}
\newtheorem{Theorem}[Lemma]{Theorem}
\newtheorem{Proposition}[Lemma]{Proposition}
\newtheorem{Corollary}[Lemma]{Corollary}

\usepackage[a4paper,top=3cm,bottom=2cm,left=3cm,right=3cm,marginparwidth=1.75cm]{geometry}

\usepackage{float}
\usepackage{nccmath}
\usepackage{mathtools}
\usepackage{amsfonts}
\usepackage{amsmath}
\usepackage{graphicx}
\usepackage[colorinlistoftodos]{todonotes}
\usepackage[colorlinks=true, allcolors=blue]{hyperref}

\title{A note on the distribution of prime ideals in real quadratic fields}



\author{Stephan~Baier}
\address{Stephan~Baier\\
	Ramakrishna Mission Vivekananda Educational Research Institute\\
	Department of Mathematics\\
	G.\ T.\ Road, PO~Belur Math, Howrah, West Bengal~711202\\
	India}
\email{stephanbaier2017@gmail.com}
\urladdr{https://www.researchgate.net/profile/Stephan\_Baier2}

\author{Sayantan Roy}
\address{Sayantan Roy\\
	Ramakrishna Mission Vidyamandira\\
	Department of Mathematics\\
	G.\ T.\ Road, PO~Belur Math, Howrah, West Bengal~711202\\
	India}
\email{roysayantan6@gmail.com}
\urladdr{}

\begin{document}
\maketitle

\begin{abstract} In this note, we give a summary of the article ``The distribution of prime ideals of imaginary quadratic fields'' by G. Harman, A. Kumchev and P. A. Lewis and establish analogous results for real quadratic fields based on the same method.
\end{abstract}

\bigskip

\tableofcontents

\section{Review of prime ideals in small regions for imaginary quadratic fields} \label{review}
The celebrated prime number theorem states that 
$$
\pi(x)\sim \frac{x}{\log x} \mbox{ as } x\rightarrow \infty, 
$$ 
where $\pi(x)$ denotes the number of primes not exceeding $x$. It is an old question if an asymptotic of this form, or at least a lower bound of expected order of magnitude, holds for primes in small intervals. In this direction, using density estimates for the zeros of the Riemann zeta function, Huxley \cite{Hux} established that 
\begin{equation} \label{asympshort}
\pi(x)-\pi(x-y) \sim \frac{y}{\log x}
\end{equation}
if $x^{7/12+\varepsilon}\le y\le x$. Using sieve methods, lower bounds of the form
\begin{equation} \label{lowershort}
\pi(x)-\pi(x-y)\gg \frac{y}{\log x}
\end{equation}
have been established for shorter intervals of length $y=x^{\theta}$ with specific exponents $\theta$ smaller than $7/12=0.5833$. The current record is due to Baker, Harman and Pintz \cite{BHP} who showed that $\theta=1/2+1/40=0.525$ is admissible. The Riemann hypothesis would imply that \eqref{asympshort} holds for $y=x^{\theta}$ with $\theta$ being any fixed real number greater than $1/2$. Substantial progress in the direction of exponents smaller than $1/2$ has not been made yet, although it may be expected that \eqref{asympshort} holds for $y=x^{\theta}$ with any fixed exponent $\theta>0$. Under the Riemann hypothesis, this is indeed true in an almost-all sense (see \cite{Hea}). Also the old conjecture that there is always a prime between any two consecutive squares is hitherto unresolved. This would be true (at least for large enough pairs of consecutive squares) if one could establish \eqref{lowershort} with $y=2\sqrt{x}$. 
  
The problem of detecting primes in short intervals has analogues for number fields $K$. The direct translation of this problem into the number field setting is to count prime ideals with norm in an interval of the form $(x-y,x]$ in the ring of integers $\mathcal{O}_K$ of $K$. Similar analytic and sieve methods as in the case $K=\mathbb{Q}$ may be applied to tackle this question. However, one may also ask more refined questions on prime ideals lying, in a sense, in small regions. To make this precise, let us assume first that $K$ has class number one. Let $\sigma_1,...,\sigma_r$ be the real and $\sigma_{r+1},\overline{\sigma_{r+1}},...,\sigma_{r+s},\overline{\sigma_{r+s}}$ be the complex embeddings, where $r+2s=n$, the degree of the field extension $K:\mathbb{Q}$. In the class number one case, every prime ideal is generated by a prime element, and we may want to count prime ideals $\mathfrak{p}$ generated by prime elements $p$ such that the $n$-tuple $\left(\sigma_1(p),...,\sigma_r(p),\Re(\sigma_{r+1}(p)),\Im(\sigma_{r+1}(p)),...,\Re(\sigma_{r+s}(p)),\Im(\sigma_{r+s}(p))\right)$ lies in a small region $\mathcal{R}\subset \mathbb{R}^n$. To tackle questions like this, one may use Hecke Gr\"o\ss encharaktere which allow to pick out ideals belonging to regions by Fourier analytic means. It is not difficult to extend the above setup to number fields with class number greater than 1. To this end, for any ideal class $\mathcal{C}$, one fixes one ideal $\mathfrak{a}_0\in \mathcal{C}^{-1}$ and counts prime ideals $\mathfrak{p}\in \mathcal{C}$ such that the principal ideal $\mathfrak{p}\mathfrak{a}_0$ is generated by an algebraic integer $a$ satisfying  
$$
\left(\sigma_1(a),...,\sigma_r(a),\Re(\sigma_{r+1}(a)),\Im(\sigma_{r+1}(a)),...,\Re(\sigma_{r+s}(a)),\Im(\sigma_{r+s}(a))\right)\in \mathcal{R}.
$$ 
General results in this direction were obtained by Coleman \cite{Col} using zero density estimates for Hecke $L$-functions.

An illuminating example is that of the number field $K=\mathbb{Q}(i)$ in which case the problem is the the same as counting prime elements $p$ (Gaussian primes) in small  regions in the complex plane (note that $\mathbb{Q}(i)$ has class number one). Most convenient is to focus on regions given by small sections of small sectors, i.e. regions of the form
$$
\mathcal{R}:=\{z\in \mathbb{C} : |z|\in I, \ \arg(z)\in J \bmod{2\pi}\},
$$
where $I$ and $J$ are short intervals. (Here the meaning of $\arg(z)\in J \bmod{2\pi}$ is that $\arg(z)-2\pi k\in J$ for a suitable $k\in \mathbb{Z}$.) Thus, one aims to detect Gaussian primes $p$ simultaneously satisfying two conditions of the form
$$
x-y\le \mathcal{N}(p)\le x \quad \mbox{and} \quad \arg(p)\in [\phi_0,\phi_0+\phi] \bmod 2\pi,
$$
where $x^{\theta_1}\le y\le x$, $\phi_0\in \mathbb{R}$ and
$x^{-\theta_2}\le \phi<2\pi$ for suitable $\theta_1,\theta_2>0$. Here $\mathcal{N}(p)$ denotes the norm and $\arg(p)$ the argument of $p$. 

Roughly speaking, the condition $x-y\le \mathcal{N}(p)\le x$ may be picked out using Mellin transform and the condition $\arg(z)\in [\phi_0,\phi_0+\phi]\bmod 2\pi$ using Hecke Gr\"o\ss encharaktere. These are characters of the multiplicative group of non-zero fractional ideals $\mathfrak{a}$ for the ring of integers $\mathcal{O}_K$. In the case $K=\mathbb{Q}(i)$, they are of the form
$$
\lambda^m(\mathfrak{a})=\left(\frac{a}{|a|}\right)^{4m}, \quad m\in \mathbb{Z},
$$
where $a$ is any generator of $\mathfrak{a}$ (see \cite[page 62, Exercise 1]{Iwa}). The right-hand side is indeed independent of the generator and hence $\lambda^m(\mathfrak{a})$ is well-defined. Here the factor $4$ in the exponent $4m$ above is of relevance: this is the number of units $1,i,-1,-i$ in $\mathcal{O}_K$. The occurrence of this factor makes it reasonable to slightly modify our above setup: Rather than detecting the condition $\arg(p)\in [\phi_0,\phi_0+\phi]\bmod 2\pi$ we detect the condition 
$\arg(p^4)\in [\phi_0,\phi_0+\phi]\bmod 2\pi$. To formulate our problem as a counting problem on prime ideals, we set
$$
\mu_\mathfrak{\mathfrak{a}}=\left(\frac{a}{|a|}\right)^4,
$$
$a$ being any generator of the ideal $\mathfrak{a}$. Now we aim to detect prime ideals $\mathfrak{p}$ simultaneously satisfying the two conditions 
\begin{equation} \label{conditions}
x-y\le \mathcal{N}(\mathfrak{p})\le x \quad \mbox{and} \quad  \arg(\mu_{\mathfrak{p}})\in [\phi_0,\phi_0+\phi] \bmod 2\pi,
\end{equation}
where $\mathcal{N}(\mathfrak{p})$ is the norm of the ideal $\mathfrak{p}$.

\begin{figure}[h]
	\includegraphics[scale=0.5]{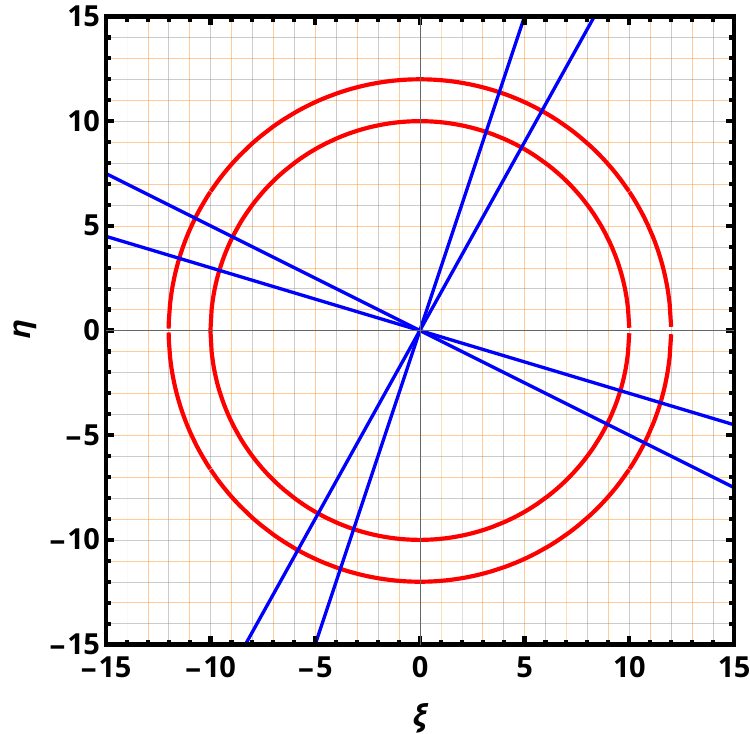}
	\caption{Conditions in (3) for $K=\mathbb{Q}(i)$}
	\label{SR_Circle}
\end{figure}

In the graphic shown in Fig[\ref{SR_Circle}], the large red circle $\xi^2+\eta^2=x$ and the small red circle $\xi^2+\eta^2=x-y$ include an annulus which contains the elements of $\mathbb{Q}(i)$ with norm in $[x-y,x]$. The blue lines include four sectors whose union contains the elements $a$ of $\mathbb{Q}(i)$ such that $(a/|a|)^4\in [\phi_0,\phi_0+\phi] \bmod 2\pi$.  

The above setup can easily be generalized to imaginary quadratic fields $K$ with class number one. In this case, the Hecke Gr\"o\ss encharaktere on the ideals $\mathfrak{a}$ for the ring of integers $\mathcal{O}_K$ are of the form
$$
\lambda^m(\mathfrak{a})=\left(\frac{a}{|a|}\right)^{wm}, \quad m\in \mathbb{Z},
$$
where $a$ is any generator of the ideal $\mathfrak{a}$ and $w$ is the number of units in $\mathcal{O}_K$, which is equal to 2,4 or 6. (For example, in the case of the field $K=\mathbb{Q}(\sqrt{-3})$, we have $w=6$.) Similarly as before, we set 
\begin{equation} \label{xidef}
\mu_{\mathfrak{a}}:=\left(\frac{a}{|a|}\right)^w
\end{equation}
and aim to detect prime ideals $\mathfrak{p}$ simultaneously satisfying the two conditions in \eqref{conditions}. 

It is easy to drop the condition of class number one. In the case of general imaginary quadratic fields, one proceeds as follows: For any ideal class $\mathcal{C}$, one fixes one ideal $\mathfrak{a}_0\in \mathcal{C}^{-1}$. 
Then the Hecke Gr\"o\ss encharaktere on the ideals $\mathfrak{a}\mathfrak{a}_0$ with $\mathfrak{a}\in \mathcal{C}$ satisfy 
$$
\lambda^m(\mathfrak{a}\mathfrak{a}_0)=\left(\frac{a}{|a|}\right)^{wm}, \quad m\in \mathbb{Z},
$$
where $a$ is any generator of the principal ideal $\mathfrak{a}\mathfrak{a}_0$ and $w$ is the number of units in $\mathcal{O}_K$, as previously. As before, we define $\mu_{\mathfrak{a}}$ by the relation \eqref{xidef} and aim to detect prime ideals $\mathfrak{p}\in \mathcal{C}$ satisfying \eqref{conditions}. 

The following is \cite[Theorem 2]{HLK}, stated in slightly modified form.

\begin{Theorem}[Harman-Lewis-Kumchev] \label{HLKresult1} Let $K$ be an imaginary quadratic number field and $\mathcal{C}$ an ideal class of $K$. Fix an ideal $\mathfrak{a}_0\in \mathcal{C}^{-1}$. Define $\theta_1:=0.765$ and $\theta_2:=1-\theta_1=0.235$. Then there is an $x_0>0$ such that if $x\ge x_0$, $x^{\theta_1}\le y\le x$, $\phi_0\in \mathbb{R}$ and $x^{-\theta_2}\le \phi< 2\pi$, then 
$$
\sum\limits_{\mathfrak{p}\in \mathcal{A}} 1 \gg \frac{\phi y}{\log x},
$$   
where 
$$
\mathcal{A}:=\{ \mathfrak{a}\in \mathcal{C} : x-y\le \mathcal{N}(\mathfrak{a})\le x \mbox{ and } \arg(\mu_{\mathfrak{a}})\in [\phi_0,\phi_0+\phi] \bmod 2\pi\}.
$$
\end{Theorem} 

Here, as throughout the sequel, we reserve the symbol $\mathfrak{p}$ for prime ideals. 
As a consequence, Harman, Lewis and Kumchev deduced the following result on representations of rational primes by indefinite quadratic forms (see \cite[Theorem 1]{HLK}).

\begin{Theorem}[Harman-Lewis-Kumchev] \label{HLKresult2} Let $\Delta<0$ be the discriminant of an imaginary quadratic field $K$, and let $\mathbb{Q}(\xi,\eta)\in \mathbb{Z}[\xi,\eta]$ be a positive definite quadratic form with discriminant $\Delta$. Then for every pair $
(s,t)\in \mathbb{R}^2$, there is another pair $(m,n)\in \mathbb{Z}^2$ for which $Q(m,n)$ is prime and 
\begin{equation} \label{sizecondition}
Q(s-m,t-n)\ll Q(s,t)^{0.53}+1.
\end{equation}
The implied constant depends only on $K$. 
\end{Theorem}

This supersedes a result by Coleman \cite{Col2} who obtained Huxley's exponent $7/12+\varepsilon$ in place of $0.53$. 
We remark that in the case of positive definite quadratic forms $Q(\xi,\eta)$ of fixed discriminant $\Delta$, the condition \eqref{sizecondition} can be reformulated as
\begin{equation*} 
||(s,t)-(m,n)||\ll ||(s,t)||^{0.53}+1,
\end{equation*}
where $||.||$ is the Euclidean norm on $\mathbb{R}^2$.

\section{Gr\"o\ss encharaktere for real quadratic fields}
Now we want to set up the problem in an analogous way for real quadratic fields $K$. Let $\sigma_1(a+b\sqrt{d})=a+b\sqrt{d}$ and $\sigma_2(a+b\sqrt{d})=a-b\sqrt{d}$ be the two embeddings of $K$ and $\epsilon>1$ be a fundamental unit in $\mathcal{O}_K$. Suppose first that $K$ has class number one. According to Hecke's work \cite{Hec1}, the Gr\"o\ss encharaktere are in this case of the form
\begin{equation*} 
\lambda^m(\mathfrak{a})=\exp\left(2\pi i m\cdot \frac{\log|\sigma_1(a)/\sigma_2(a)|}{2\log \epsilon}\right), \quad m\in \mathbb{Z},
\end{equation*}
where $a\in \mathcal{O}_K$ is any generator of the (principal) ideal $\mathfrak{a}$. Again, it is easy to check that the right-hand side only depends on $\mathfrak{a}$ and hence $\lambda^m(\mathfrak{a})$ is well-defined. Obviously, $\lambda^m(\mathfrak{a})=\mu_{\mathfrak{a}}^m$, where
\begin{equation} \label{etadef}
\mu_{\mathfrak{a}}:=\exp\left(2\pi i\cdot \frac{\log|\sigma_1(a)/\sigma_2(a)|}{2\log \epsilon}\right). 
\end{equation}
This takes now the role of $\mu_{\mathfrak{a}}$ in the case of imaginary quadratic fields. 
So similarly as before, our problem becomes to detect prime ideals $\mathfrak{p}$ satisfying conditions as in 
\eqref{conditions}. 

The necessary modifications in the case of a class number greater than one are carried out similarly as for imaginary quadratic fields: For any ideal class $\mathcal{C}$, one fixes one ideal $\mathfrak{a}_0\in \mathcal{C}^{-1}$. 
Then the Hecke Gr\"o\ss encharaktere on the ideals $\mathfrak{a}\mathfrak{a}_0$ with $\mathfrak{a}\in \mathcal{C}$ satisfy
\begin{equation} \label{lambdamrel}
\lambda^m(\mathfrak{a}\mathfrak{a}_0)=\mu_{\mathfrak{a}}^m,
\end{equation} 
where $\mu_{\mathfrak{a}}$ is defined as in \eqref{etadef} and 
$a$ is any generator of the principal ideal $\mathfrak{a}\mathfrak{a}_0$.

\section{Main results}
Our first main result in this article is that Theorem \ref{HLKresult1} remains valid in the context of real quadratic fields. 

\begin{Theorem} \label{mainresult1} Theorem \ref{HLKresult1} above holds for real quadratic fields $K$, where our notations are as in the previous section.  
\end{Theorem} 

As a consequence, we will deduce the following analogue of Theorem \ref{HLKresult2} on representations of rational primes by definite quadratic forms. 

\begin{Theorem} \label{mainresult2}  Assume that $\delta>0$. Let $\Delta$ be the discriminant of $K$ and $\mathbb{Q}(\xi,\eta)=a\xi^2+b\xi\eta+c\eta^2\in \mathbb{Z}[\xi,\eta]$ be an indefinite quadratic form with discriminant $\Delta$. Let $(s,t)\in \mathbb{R}^2$ such that $Q(s,t)>0$ and the angles between each of the asymptotes of the hyperbola $Q(\xi,\eta)=1$ and the line passing through the origin and the point $(s,t)$ are greater than or equal to $\delta$. Then there is another pair $(m,n)\in \mathbb{Z}^2$ for which $Q(m,n)$ is prime and 
\begin{equation*} \label{sizeconditionreformulated}
||(s,t)-(m,n)||\ll  ||(s,t)||^{0.53}+1,
\end{equation*}
where $||.||$ is the Euclidean norm on $\mathbb{R}^2$. The implied constant depends only on the quadratic form and $\delta$. 
\end{Theorem}
$ $\\
{\bf Remark:} The asymptotes of the hyperbola $Q(s,t)=1$ (and all other hyperbolas $Q(\xi,\eta)=N>0$) are the lines passing through the origin and the points $(-b\pm \sqrt{\Delta},2a)$, respectively.\\ \\
{\bf Acknowledgements:} The authors would like to thank the Ramakrishna Mission Vivekananda Educational and Research Institute for providing excellent working conditions. They would also like to thank the anonymous referee for valuable comments. 
  
\section{Description of the Harman-Lewis-Kumchev method}
In the following, we describe the method used by Harman, Lewis and Kumchev \cite{HLK} to prove Theorem \ref{HLKresult1}. This method will also yield our first main result, Theorem \ref{mainresult1}, as we will explain later. The aim is to prove that 
\begin{equation} \label{maingoal}
\sum\limits_{\mathfrak{p}} I_{\mathcal{A}}(\mathfrak{p})\gg \frac{\phi y}{\log x},
\end{equation}
where $I_{\mathcal{A}}$ is the indicator function of the set
$$
\mathcal{A}:=\{ \mathfrak{a}\in \mathcal{C}: x-y\le \mathcal{N}(\mathfrak{a})\le x \mbox{ and } \arg(\mu_{\mathfrak{a}})\in [\phi_0,\phi_0+\phi] \bmod 2\pi \}.
$$
Here we note that the volume of the set 
$$
\{(t,\theta): x-y\le t\le x, \ \phi_0\le \theta\le \phi_0+\phi\} \subset \mathbb{R}^2
$$
equals $\phi y$. Since the probability of a randomly chosen ideal $\mathfrak{a}$ to be prime is of size about $1/\log \mathcal{N}(\mathfrak{a})$ by Landau's prime ideal theorem, and the prime ideals are uniformly distributed over the classes $\mathcal{C}$, we expect approximately 
$$
\frac{\phi y}{h \log x}
$$ 
prime ideals in $\mathcal{A}$, where $h$ is the class number of $K$.  

To make $I_{\mathcal{A}}$ more suitable for an application of Fourier analysis, it is useful to smooth this function out. Harman, Lewis and Kumchev \cite{HLK} approximated it by a function 
\begin{equation} \label{Psidefi}
\Psi(\mathfrak{a}):=\psi_1(\mathcal{N}(\mathfrak{a}))\psi_2(\arg(\mu_{\mathfrak{a}})),
\end{equation}
where $\psi_1$ and $\psi_2$ satisfy the following conditions. 

Let $\Delta_1:=yx^{-\eta}$, where $\eta>0$ is small enough. The function $\psi_1(t)$ is supposed to be a function of the $C^{\infty}(\mathbb{R})$ class such that
\begin{itemize}
\item $\psi_1(t)=1$ if $x-y+\Delta_1\le t\le x-\Delta_1$,
\item $\psi_1(t)=0$ if $t\not\in (x-y,x)$,
\item $0\le \psi_1(t)\le 1$ if $x-y\le t\le x$,
\item $\psi_1^{(j)}\ll j^{2j}\Delta_1^{-j}$ for $j=1,2,...$.
\end{itemize}
Let $\Delta_2:=\phi x^{-\eta}$ and $r=[2/\eta]+1$. The function $\psi_2(t)$ is supposed to be a $2\pi$-periodic function such that
\begin{itemize}
\item $\psi_2(t)=1$ if $\phi_0+\Delta_2\le t\le \phi_0+\phi-\Delta_2$,
\item $\psi_2(t)=0$ if $\phi_0+\phi\le t\le \phi_0+2\pi$,
\item $\psi_2(t)=0$ if $\phi_0\le t\le \phi_0+\phi$,
\item $\psi_2(t)$ has a Fourier expansion 
\begin{equation} \label{Fourier}
\psi_2(t)=\hat{\psi}_2(0)+\sum\limits_{m\not=0} \hat{\psi}_2(m) e^{imt},
\end{equation}
where $\hat{\psi}_2(0)=(\phi-\Delta_2)/2\pi$ and, for $m\not=0$,
\begin{equation} \label{Fourierbound}
|\hat{\psi}_2(m)|\le \min\left(\phi,\frac{2}{\pi |m|}, \frac{2}{\pi m}\left(\frac{2r}{|m|\Delta_2}\right)^r\right). 
\end{equation}
\end{itemize}

Clearly, $I_\mathcal{A}(\mathfrak{a})\ge \Psi(\mathfrak{a})$ for all ideals $\mathfrak{a}$. Thus it suffices to prove that
$$
\sum\limits_{\mathfrak{p}} \Psi(\mathfrak{p})\gg \frac{\phi y}{\log x}. 
$$
The sieve method used by Harman, Lewis and Kumchev follows the general strategy of Harman's sieve to compare the number of prime ideals in $\mathcal{A}$ with the number of prime ideals in a much larger set $\mathcal{B}$ for which the set of prime ideals contained in it is under good control. The setup is somewhat similar to that in the work of Baker, Harman and Pintz \cite{BHP} on primes in short intervals. Here the choice of $\mathcal{B}$ is 
\begin{equation} \label{Bdef}
\mathcal{B}=\{\mathfrak{a}\in C : x-y_1\le \mathcal{N}(\mathfrak{a})\le x\},
\end{equation}
where $y_1:=x\exp\left(-3\log^{1/3} x\right)$. (In the case of $K=\mathbb{Q}(i)$, this corresponds to the set of prime elements in a large annulus.) Using Hecke $L$-functions with class group characters to pick out the condition $\mathfrak{a}\in \mathcal{C}$, it is easy to establish that
\begin{equation} \label{Bprimes}
\sum\limits_{\mathfrak{p}\in \mathcal{B}} 1\sim \frac{y_1}{h\log x} 
\end{equation}
by familiar analytic techniques, making use of zero-free regions for these $L$-functions. Here we note that $y_1$ is sufficiently large for these techniques to work. 

Following the general idea of Harman's sieve, one compares the quantities
$$
S(\mathcal{A},z):=\sum\limits_{\substack{\mathfrak{a}\in \mathcal{A}\\ (\mathfrak{a},\mathfrak{P}(x))=1}} \Psi(\mathfrak{a}) \quad \mbox{and} \quad S(\mathcal{B},z):=\sum\limits_{\substack{\mathfrak{a}\in \mathcal{B},\\ (\mathfrak{a},\mathfrak{P}(x))=1}} 1,
$$
where $z\ge 2$ and 
$$
\mathfrak{P}(z):=\prod\limits_{\mathcal{N}(\mathfrak{p})<z} \mathfrak{p}.
$$
Clearly, 
$$
S\left(\mathcal{A},x^{1/2}\right)=\sum\limits_{\mathfrak{p}\in \mathcal{A}} \Psi(\mathfrak{a})+O\left(x^{1/2}\right) 
$$
and 
\begin{equation} \label{Bsqareroot} 
S\left(\mathcal{B},x^{1/2}\right)=\sum\limits_{\mathfrak{p}\in \mathcal{B}} 1+O\left(x^{1/2}\right).
\end{equation}
By a number of Buchstab decompositions, the quantity $S\left(\mathcal{A},x^{1/2}\right)$ is now expressed in the form
\begin{equation} \label{Adecomp}
S\left(\mathcal{A},x^{1/2}\right)=\sum\limits_{j=1}^k S_j - \sum\limits_{j=k+1}^l S_j,
\end{equation}
where $S_j\ge 0$ are certain sifting sums. Analogous Buchstab decompositions lead to an expression of the form
\begin{equation} \label{Bdecomp}
S\left(\mathcal{B},x^{1/2}\right)=\sum\limits_{j=1}^k \tilde{S}_j - \sum\limits_{j=k+1}^l \tilde{S}_j,
\end{equation}
where $\tilde{S}_j\ge 0$. Now the right-hand sides of \eqref{Adecomp} and \eqref{Bdecomp} are being compared. It turns out that for a suitable $k'<k$, asymptotic formulas of the form
\begin{equation} \label{asympAB}
S_j=\frac{\phi y}{2\pi y_1} \tilde{S}_j (1+o(1))
\end{equation}
can be established if $1\le j\le k'<k$ or $j>k$. Note that the factor $(\phi y)/(2\pi y_1)$ on the right-hand side equals the ratio of the volumes of the sets $[x-y,x]\times [\phi_0,\phi_0+\phi]$ and $[x-y_1,x]\times (0,2\pi]$ in $\mathbb{R}^2$. The terms $S_j$ with $k'<j\le k$ are discarded from the right-hand side of \eqref{Adecomp}, giving an inequality of the form
\begin{equation} \label{finalin}
\begin{split}
S\left(\mathcal{A},x^{1/2}\right)\ge\sum\limits_{j=1}^{k'} S_j - \sum\limits_{j=k+1}^l S_j = & \frac{\phi y}{2\pi y_1}\left(\sum\limits_{j=1}^{k'} \tilde{S}_j - \sum\limits_{j=k+1}^l \tilde{S}_j \right)(1+o(1))\\ = & 
\frac{\phi y}{2\pi y_1}\left(S\left(\mathcal{B},x^{1/2}\right)-\sum\limits_{j=k'+1}^k \tilde{S}_j\right)(1+o(1)).
\end{split}
\end{equation}
It then remains to be shown that the last line in \eqref{finalin} is actually positive and of expected order of magnitude $\gg \phi y/\log x$. To this end, one bounds $\tilde{S}_j$ for $j\in \{k'+1,...,k\}$ from above, which requires careful numerical calculations, and uses \eqref{Bprimes} and \eqref{Bsqareroot} to approximate $S\left(\mathcal{B},x^{1/2}\right)$.

To establish asymptotic formulas as in \eqref{asympAB}, one needs to connect sums of the form $\sum\limits_{\mathfrak{a}\in \mathcal{A}} c(\mathfrak{a})\Psi(\mathfrak{a})$ with corresponding sums of the form 
$\sum\limits_{\mathfrak{a}\in \mathcal{B}} c(\mathfrak{a})$ to obtain 
asymptotic formulas of the form
\begin{equation} \label{asy1}
\sum\limits_{\mathfrak{a}\in \mathcal{A}} c(\mathfrak{a})\Psi(\mathfrak{a})=
\frac{\phi y}{2\pi y_1}(1+o(1))
\sum\limits_{\mathfrak{a}\in \mathcal{B}} c(\mathfrak{a}).
\end{equation}
Here $c(\mathfrak{a})$ are specific divisor bounded coefficients, meaning that 
$c(\mathfrak{a})\ll \tau(\mathfrak{a})^D$ for some power $D>0$, where $\tau(\mathfrak{a})$ is the number integral ideal divisors of $\mathfrak{a}$. At some point in the method, this will not be quite sufficient, and one needs to establish similar asymptotic formulas for bilinear sums with one short range of the form
\begin{equation} \label{asy2}
\sum\limits_{\substack{\mathfrak{a},\mathfrak{b}\\ \mathfrak{ab}\in \mathcal{A}\\ \mathcal{N}(\mathfrak{a})\le x^{\eta}}} c(\mathfrak{a})\Psi(\mathfrak{ab})=
\frac{\phi y}{2\pi y_1}(1+o(1))
\sum\limits_{\substack{\mathfrak{a},\mathfrak{b}\\ \mathfrak{ab}\in \mathcal{B}\\ \mathcal{N}(\mathfrak{a})\le x^{\eta}}} c(\mathfrak{a}). 
\end{equation}

The method of Harman, Lewis and Kumchev first reduces the said sums to Dirichlet polynomials with Gr\"o\ss encharaktere. Here the Mellin transform is applied to treat the function $\psi_1$ and the Fourier series expansion in \eqref{Fourier} to treat $\psi_2$. (Recall that $\psi_1$ and $\psi_2$ are implicit in the definition of $\Psi$ in \eqref{Psidefi}.) This Fourier series expansion can conveniently be cut off using the bound \eqref{Fourierbound} for the Fourier coefficients. Now asymptotic formulas as in \eqref{asy1} and \eqref{asy2} are achieved by using mean value theorems for the said Dirichlet polynomials. These mean value estimates rely on results of Coleman in \cite{Col3}, \cite{Col2} and \cite{Col} and results of Harman, Baker and  Pintz in \cite{BH} and \cite{BHP}. 

A careful examination reveals that both the sieve method and the estimates for Dirichlet polynomials work in the same way for real quadratic as for imaginary quadratic fields. Indeed, all underlying lemmas apply to general number fields. Also the analytic treatments of $\psi_1$ and $\psi_2$ remain unaltered. Therefore, the method of Harman, Lewis and Kumchev gives our first main result, Theorem \ref{mainresult1}, in exactly the same way as their corresponding Theorem \ref{HLKresult1}. The only difference is that $\mu_\mathfrak{a}$ is in the real quadratic case defined as in \eqref{etadef}, which differs from its definition \eqref{xidef} in the imaginary quadratic case. In the next section, we will describe how the sums $\sum\limits_{\mathfrak{a}\in \mathcal{A}} c(\mathfrak{a})\Psi(\mathfrak{a})$ and $\sum\limits_{\substack{\mathfrak{ab}\in \mathcal{A}\\ \mathcal{N}(\mathfrak{a})\le x^{\eta}}} c(\mathfrak{a})\Psi(\mathfrak{ab})
$ (and the corresponding $\mathcal{B}$-sums) are reduced to Dirichlet polynomials and which mean value estimates for Dirichlet polynomials are needed. This is a sketch of \cite[section 3]{HLK}. We will not go into the rather complicated details of \cite[sections 4-6]{HLK} which contain the proofs of the required mean value estimates and the sieve method. All what we point out is that all details in these sections apply to real quadratic fields as well.  

In the last section of this paper, we will deduce our Theorem \ref{mainresult2} on rational primes represented by indefinite quadratic forms with variables restricted to short intervals from our Theorem \ref{mainresult1}. 

\section{Reduction to mean value estimates for Dirichlet polynomials}
Employing the Fourier series expansion of $\psi_2$ in \eqref{Fourier}, equation \eqref{lambdamrel} and $|\mu_{\mathfrak{a}}|=1$, we have 
$$
\psi_2(\arg(\mu_{\mathfrak{a}}))=\sum\limits_{m\in \mathbb{Z}} \hat{\psi}_2(m) e^{im \arg(\mu_{\mathfrak{a}})}=\sum\limits_{m\in \mathbb{Z}} \hat{\psi}_2(m) \lambda^m(\mathfrak{a}\mathfrak{a}_0).  
$$
Using \eqref{Fourierbound}, the right-hand side can be truncated at $M:=\left[\Delta_2^{-1}x^{\eta}\right]+1$ at the cost of an error of size $O(x^{-2})$, i.e.,
$$
\psi_2(\arg(\mu_{\mathfrak{a}}))=\sum\limits_{|m|\le M} \hat{\psi}_2(m) \lambda^m(\mathfrak{a}\mathfrak{a}_0)+O\left(x^{-2}\right). 
$$ 
By the definition of $\Psi(\mathfrak{a})$ in \eqref{Psidefi} and multiplicativity of $\lambda^m$, we can therefore write
$$
\sum\limits_{\mathfrak{a}\in \mathcal{A}} c(\mathfrak{a})\Psi(\mathfrak{a})=
\sum\limits_{\mathfrak{a}\in \mathcal{C}} c(\mathfrak{a})\Psi(\mathfrak{a})=
\sum\limits_{|m|\le M} \hat{\psi}_2(m) \lambda^m(\mathfrak{a}_0) \sum\limits_{\mathfrak{a}\in \mathcal{C}} c(\mathfrak{a})\lambda^m(\mathfrak{a})\psi_1(\mathcal{N}(\mathfrak{a}))+O(1). 
$$
To handle the function $\psi_1$, we use its Mellin transform
$$
\hat{\psi}_1(s)=\int\limits_0^{\infty} \psi_1(t)t^{s-1}dt.
$$
By the Mellin inversion formula, we may express $\psi_1(t)$ in terms of $\hat{\psi}_1(s)$ in the form
$$
\psi_1(t)=\frac{1}{2\pi i} \int\limits_{1/2-i\infty}^{1/2+i\infty} \hat{\psi}_1(s)t^{-s}ds. 
$$
It follows that
\begin{equation} \label{afterMellin}
\sum\limits_{\mathfrak{a}\in \mathcal{A}} c(\mathfrak{a})\Psi(\mathfrak{a})=
\frac{1}{2\pi i} \sum\limits_{|m|\le M} \hat{\psi}_2(m) \lambda^m(\mathfrak{a}_0) \int\limits_{1/2-i\infty}^{1/2+i\infty} F(s,\lambda^m)\hat{\psi}_1(s)ds+O(1),
\end{equation}
where $F(s,\lambda^m)$ is the Dirichlet polynomial
$$
F(s,\lambda^m):=\sum\limits_{\mathfrak{a}\in C} c(\mathfrak{a})\lambda^m(\mathfrak{a})\mathcal{N}(\mathfrak{a})^{-s}. 
$$
Here the sum on the right-hand side can be restricted to $\mathfrak{a}\in \mathcal{C}$ with $\mathcal{N}(\mathfrak{a})$ lying in an interval of the form $[c_1 x,c_2x]$, where $c_1,c_2$ are suitable constants satisfying $c_1<1<c_2$. (Note that $\psi_1(t)$ equals zero outside the interval $[x-y,x]$.) Moreover, bounding $\hat{\psi}_1(s)$ appropriately (for the details, see \cite[section 3]{HLK}), the integral on the right-hand side of \eqref{afterMellin} can be truncated at $\Im s=T_1:=x^{1+\eta}\Delta_1^{-1}$ at the cost of an error of size $O(1)$, i.e.,
\begin{equation} \label{afterMellinTruncation}
\sum\limits_{\mathfrak{a}\in \mathcal{A}} c(\mathfrak{a})\Psi(\mathfrak{a})=
\frac{1}{2\pi i} \sum\limits_{|m|\le M} \hat{\psi}_2(m) \lambda^m(\mathfrak{a}_0) \int\limits_{1/2-iT_1}^{1/2+iT_1} F(s,\lambda^m)\hat{\psi}_1(s)ds+O(1).
\end{equation}
Now one isolates the main term contribution 
\begin{equation} \label{mainterm} 
{\bf M}_{\mathcal{A}}=\frac{\hat{\psi}_2(0)}{2\pi i} \int\limits_{1/2-iT_1}^{1/2+iT_1} F(s,\lambda^0)\hat{\psi}_1(s)ds
\end{equation}
coming from $m=0$ and bounds the error by 
$$
\ll A_1A_2E+O(1),
$$
where 
$$
A_1:=\max\limits_{0< |m|\le M} |\hat{\psi}_1(m)|, \quad A_2:=
\sup\limits_{|t|\le T_1} \left|\hat{\psi}_2\left(\frac{1}{2}+it\right)\right|
$$
and 
\begin{equation} \label{Edef}
E:=\sum\limits_{0<|m|\le M} \int\limits_{-T_1}^{T_1} \left|F\left(\frac{1}{2}+it,\lambda^m\right)\right|dt.
\end{equation}
The error $A_1A_2E$ turns out to be smaller than the main term ${\bf M}_{\mathcal{A}}$ if 
$E$ is slightly smaller than $x^{1/2}$. Indeed, a mean value estimate for Dirichlet polynomials of the form
\begin{equation} \label{Dirineeded1}
E\ll x^{1/2}\exp\left(-B(\log x)^{1/4}\right)
\end{equation}
for a suitable constant $B>0$ does the job (see \cite[section 3]{HLK}). Bounding $A_1$ and $A_2$ suitably (for the details see \cite[section 3]{HLK}) then yields
\begin{equation} \label{Aapprox}
\sum\limits_{\mathfrak{a}\in \mathcal{A}} c(\mathfrak{a})\Psi(\mathfrak{a})={\bf M}_{\mathcal{A}}+O\left(\phi y \exp\left(-B(\log x)^{1/4}\right)\right). 
\end{equation}
Moreover, it turns out to be convenient to truncate the integral on the right-hand side of \eqref{mainterm} at $\Im s=\pm T_0$ with $T_0:=\exp\left(\log^{1/3} x\right)$. To this end, another mean value estimate for a Dirichlet polynomial of the form
\begin{equation} \label{Dirineeded2}
\int\limits_{T_0}^{T_1} \left|F\left(\frac{1}{2}+it,\lambda^0\right)\right|dt\ll
x^{1/2}\exp\left(-B(\log x)^{1/4}\right)
\end{equation}
is required. If this holds, then the assumption $\hat{\psi}_2(0)=(\phi-\Delta_2)/(2\pi)$ and a suitable approximation of $\hat{\psi}_1(s)$ (see \cite[section 3]{HLK}) yield 
\begin{equation} \label{mainterm'} 
{\bf M}_{\mathcal{A}}=\frac{\phi y}{(2\pi)^2 i} \int\limits_{1/2-iT_0}^{1/2+iT_0} F(s,\lambda^0)x^{s-1}ds+O\left(\phi y \exp\left(-B(\log x)^{1/4}\right)\right).
\end{equation}

Recall the definition of the set $\mathcal{B}$ in \eqref{Bdef}. A similar (but simpler) calculation leads to
\begin{equation} \label{Bapprox}
\sum\limits_{\mathfrak{a}\in \mathcal{B}} c(\mathfrak{a})={\bf M}_{\mathcal{B}}+O\left(y_1 \exp\left(-B(\log x)^{1/4}\right)\right),
\end{equation}
where 
\begin{equation} \label{Bmainterm}
{\bf M}_{\mathcal{B}}=\frac{y_1}{2\pi i} \int\limits_{1/2-iT_0}^{1/2+iT_0} F(s,\lambda^0)x^{s-1}ds +O\left(y_1 \exp\left(-B(\log x)^{1/4}\right)\right).
\end{equation}
Now comparing \eqref{Aapprox} and \eqref{Bapprox} and using \eqref{mainterm'} and \eqref{Bmainterm}, we obtain the desired relation \eqref{asy1}. 
Similar considerations establish also the relation \eqref{asy2}, where one requires the above mean value estimate \eqref{Dirineeded1} to hold for $E$ as defined in \eqref{Edef} with the Dirichlet polynomial replaced by
$$
F(s,\lambda^m):=\sum\limits_{\substack{\mathfrak{a},\mathfrak{b}\\ \mathfrak{ab}\in \mathcal{C}\\ \mathcal{N}(\mathfrak{a})\le x^{\eta}\\ c_1x\le \mathcal{N}(\mathfrak{ab})\le c_2x}} c(\mathfrak{a})\lambda^m(\mathfrak{ab})\mathcal{N}(\mathfrak{ab})^{-s}. 
$$
Theorem \ref{HLKresult1} of Harman, Lewis and Kumchev then follows as described in the previous section. 

Now it remains to establish the mean value estimates \eqref{Dirineeded1} and \eqref{Dirineeded2} for the Dirichlet polynomials containing the specific coefficients $c(\mathfrak{a})$ arising from the sieve method. Both the sieve part and the mean value estimates require a large amount of calculations and are carried out in detail in \cite[sections 4-6]{HLK}.  As pointed out in the previous section of this article, our Theorem \ref{mainresult1} is established in exactly the same way, and we have therefore completed its proof. 

\section{Primes represented by quadratic forms}
\subsection{Link to quadratic fields}
Denote the set of rational primes by $\mathbb{P}$. Let a positive definite quadratic form $Q(\xi,\eta)\in \mathbb{Z}[\xi,\eta]$ be given by the equation
$$
Q(\xi,\eta)=a\xi^2+b\xi\eta+c\eta^2.
$$
This form has discriminant 
$$
\Delta=b^2-4ac<0.
$$
Let $K=\mathbb{Q}(\sqrt{\Delta})$ and 
$$
\mathfrak{d}=\left[a,\frac{b-\sqrt{\Delta}}{2}\right]
$$
be the ideal generated by $a$ and $(b-\sqrt{\Delta})/2$ in $\mathcal{O}_K$. Let further $\mathcal{C}$ be the ideal class such that $\mathfrak{d}\in \mathcal{C}^{-1}$. 
The link between prime ideals $\mathfrak{p}\in \mathcal{C}$ with norm $\mathcal{N}(\mathfrak{p})=p\in \mathbb{P}$ and representations of primes $p$ by $Q(\xi,\eta)$ was provided by the following result in the paper \cite{Col2} by Coleman. 

\begin{Proposition}[Theorem 2.2 in \cite{Col2}] If $Q(\xi,\eta)$ is positive definite, then there is a one-to-one correspondence between prime ideals $\mathfrak{p}\in \mathcal{C}$
with norm $\mathcal{N}(\mathfrak{p})=p\in \mathbb{P}$ and pairs of rational integers gcd$(m,n)$ where
$(m,n)=1$, with $m, n>0$ if $\Delta=-4$ or $-3$; and $m>0$ or $m = 0$, $n = 1$ if $\Delta<-4$. If $\mathfrak{p}\in \mathcal{C}$ and $(m,n)\in \mathbb{Z}^2$ are two correspondents, then we have
\begin{equation} \label{pd}
\mathfrak{pd}=\left(ma+n\frac{b-\sqrt{\Delta}}{2}\right)
\end{equation}
and
\begin{equation} \label{Qmn}
Q(m, n) = \mathcal{N}(\mathfrak{p})= p. 
\end{equation}
\end{Proposition}

Coleman \cite{Col2} also made the above correspondence quantitative. He demonstrated that the conditions in \eqref{conditions} translate into the conditions
$$
x-y\le |\{m,n\}|^2\le x \quad \mbox{and} \quad w\arg(\{m,n\})\in [\varphi_0,\varphi_0+\varphi] \bmod{2\pi}
$$ 
on the pair $(m,n)$, where 
$$
\{m,n\}=ma^{1/2}+n\cdot \frac{b-\sqrt{\Delta}}{2a^{1/2}}.
$$
This allows to deduce Theorem \ref{HLKresult2} from Theorem \ref{HLKresult1}. 

In the case of indefinite forms with discriminant $\Delta=b^2-4ac>0$, the above theorem remains valid except that now every prime ideal $\mathfrak{p}$ corresponds to an infinite set of pairs $(m,n)$ of integers. We have the following.

\begin{Proposition} If $Q(\xi,\eta)$ is indefinite, then there is a one-to-one correspondence between prime ideals $\mathfrak{p} \in \mathcal{C}$
with norm $\mathcal{N}(\mathfrak{p})=p\in \mathbb{P}$ and disjoint sets $\mathcal{M}$ of pairs of rational integers $(m,n)$ with gcd$(m,n)=1$ as follows. A pair of integers $(m,n)$ belongs to the set $\mathcal{M}$ corresponding to a prime ideal $\mathfrak{p}\in \mathcal{C}$ with norm $\mathcal{N}(\mathfrak{p})=p\in \mathbb{P}$ if and only if \eqref{pd} and \eqref{Qmn} hold. 
\end{Proposition}

In the following, we make this correspondence quantitative for the simplest case when $K=\mathbb{Q}(\sqrt{d})$ is a field with class number 1, where $d$ is a square-free natural number satisfying $d\equiv 2,3 \bmod{4}$. We will later indicate which modifications need to be made to cover the general case. In the above special case, the discriminant of $K$ equals $\Delta=4d$ and the quadratic form in question is equivalent to the form $\xi^2-d\eta^2$. Therefore it suffices to consider this quadratic form $Q(\xi,\eta)=\xi^2-d\eta^2$. In the following, all $\ll$-constants are allowed to depend on $d$, but we will not indicate this dependency. All other dependencies will be indicated. 

Since $\mathcal{N}(\mathfrak{p})=p=Q(m,n)$, we have 
\begin{equation} \label{firstcondi}
x-y\le \mathcal{N}(\mathfrak{p})\le x \Longleftrightarrow x-y\le Q(m,n)\le x.
\end{equation}
Further, set 
$$
f(\xi,\eta):=\xi-\eta\sqrt{d} \quad \mbox{and} \quad 
F(\xi,\eta):=\frac{\log\left|\sigma_1\left(f(\xi,\eta)\right)/\sigma_2\left(f(\xi,\eta)\right)\right|}{2\log \epsilon} \quad \mbox{for } (\xi,\eta)\in \mathbb{R}^2.
$$
Then, noting that $\mathfrak{d}=(1)$, \eqref{pd} takes the form 
\begin{equation} \label{pd2}
\mathfrak{p}=(f(m,n)). 
\end{equation}
Moreover, for any generator $\pi_0$ of $\mathfrak{p}$, there is a pair $(m,n)\in \mathcal{M}$ such that 
$$
\pi_0=f(m,n)
$$
(take a fixed generator $\pi_1=f(m_1,n_1)$ and multiply by a suitable unit) and, consequently,
$$
\sigma_i(\pi_0)=\sigma_i(f(m,n)) \quad \mbox{for } i=1,2.
$$
Taking $a=\pi_0=\pi_1\epsilon^l$ in \eqref{etadef} where $\pi_1$ is a fixed generator of $\mathfrak{p}$ and $l$ runs over $\mathbb{Z}$, it follows that $\arg(\mu_{\mathfrak{p}})\in [\phi_0,\phi_0+\phi] \bmod{2\pi}$ if and only if for some $k\in \mathbb{Z}$ there exists a pair $(m,n)$ belonging to the set $\mathcal{M}$ corresponding to $\mathfrak{p}$ such that 
$$ 
F(m,n) \in [\phi_0+2\pi k,\phi_0+\phi+2\pi k].
$$
Set
\begin{equation} \label{Sxyov}
\mathcal{S}(x,y,\omega,\varphi):= \{(\xi,\eta)\in \mathbb{R}^2\setminus \{(0,0)\} : 
x-y\le Q(\xi,\eta)\le x, \ \omega\le F(\xi,\eta) \le \omega+\phi\}.
\end{equation}
By the above considerations, Theorem \ref{mainresult1} implies the following. 

\begin{Corollary} \label{coro} Let  $x\ge 2$ large enough and $\omega\in \mathbb{R}$. Then under the conditions $x^{\theta_1}\le y\le x$ and $x^{\theta_2}\le \phi<2\pi$ in Theorems \ref{HLKresult1} and \ref{mainresult1}, there exists 
$(m,n)\in \mathcal{S}(x,y,\omega,\phi)\cap \mathbb{Z}^2$ such that $Q(m,n)=p$ is a prime. 
\end{Corollary}

\subsection{Geometric interpretation} 
Geometrically, the set $\mathcal{Q}$ of points $(\xi,\eta)$ satisfying $x-y\le Q(\xi,\eta)\le x$ is included by the two hyperbolas $Q(\xi,\eta)=x-y$ and $Q(\xi,\eta)=x$. The asymptotes of these hyperbolas are the two lines $l_{1,2}$ with slopes 
$$
\frac{\eta}{\xi}=\pm \frac{1}{\sqrt{d}}.
$$
Since 
$$
F(\xi,\eta)=F\left(1,\frac{\eta}{\xi}\right),
$$ 
the condition $F(\xi,\eta)\in [\omega,\omega+\phi]$ is equivalent to $\eta/\xi$ lying in a union of two intervals.  Hence, the set $\mathcal{F}$ of points $(\xi,\eta)$ satisfying $F(\xi,\eta)\in [\omega,\omega+\phi]$ is a union of two sectors, each of them included by two lines. Indeed, we calculate that
\begin{equation}
\begin{split}
& \omega\le \frac{\log\left|\sigma_1\left(f(\xi,\eta)\right)/\sigma_2\left(f(\xi,\eta)\right)\right|}{2\log \epsilon}\le \omega+\phi\\ \Longleftrightarrow\ 
& \frac{1}{\sqrt{d}}\cdot \tanh \left(\omega\log \epsilon\right)\le -\frac{\eta}{\xi}\le \frac{1}{\sqrt{d}}\cdot \tanh \left((\omega+\phi)\log \epsilon\right) \mbox{ or } \\
& 
\frac{1}{\sqrt{d}}\cdot \coth \left(\omega\log \epsilon\right)\le -\frac{\eta}{\xi}\le \frac{1}{\sqrt{d}}\cdot \coth \left((\omega+\phi)\log \epsilon\right).
\end{split}
\end{equation}
The sector corresponding to the first inequality for $\eta/\xi$ intersects with the above hyperbolas, the other sector does not. The intersection $\mathcal{Q}\cap \mathcal{F}$ equals the set $\mathcal{S}(x,y,\omega,\phi)$. It consists of two antipodal connected components $\mathcal{S}_1(x,y,\omega,\phi)$ and $\mathcal{S}_2(x,y,\omega,\phi)$.

\begin{figure}[h]
	\includegraphics[scale=0.6]{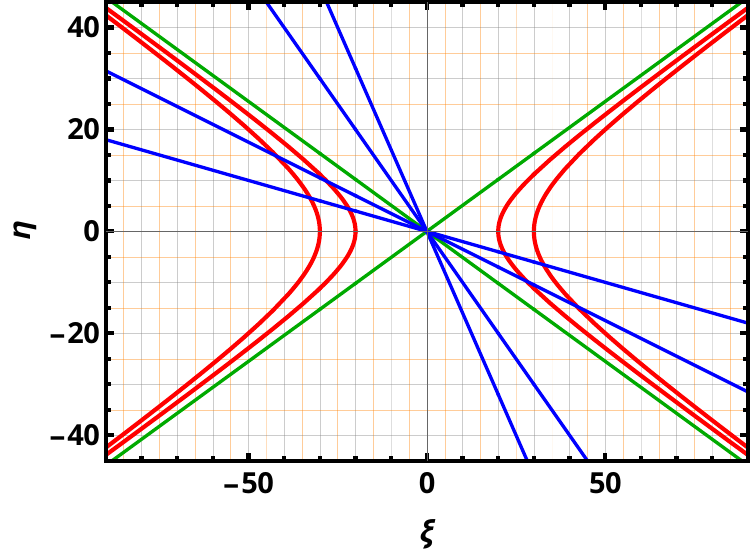}
	\caption{Sectors and hyperbolas for $Q(\xi,\eta)=\xi^2-d\eta^2$}
	\label{SR_Parabola}
	\end{figure}

Fig[\ref{SR_Parabola}] represents a diagram for the case of the quadratic form $Q(\xi,\eta)=\xi^2-2\eta^2$. 
In the graphic, the hyperbola $Q(\xi,\eta)=x$ and $Q(\xi,\eta)=x-y$ are marked red, the sectors are included by the blue lines and the asymptotes are the green lines.  

\subsection{Proof of Theorem \ref{mainresult2} for $Q(\xi,\eta)=\xi^2-d\eta^2$} Let $\delta>0$ and assume that $x$ is large enough so that Corollary \ref{coro} is applicable. Suppose that $y$ and $\phi$ satisfy the conditions in Corollary \ref{coro}. Assume $S$ is the sector included by the two asymptotes $l_1$ and $l_2$ and containing the points $(\xi,\eta)$ for which $Q(\xi,\eta)$ is non-negative. In particular, this sector contains the hyperbolas $Q(\xi,\eta)=x-y$ and $Q(\xi,\eta)=x$. Without loss of generality, assume that $\delta$ is smaller than half of the angle between $l_1$ and $l_2$. Assume further $S'$ is the closed subsector of $S$ included by two lines $l_1'$ and $l_2'$, where the angle between $l_i$ and $l_i'$ equals $\delta$ for $i=1,2$. 
Suppose that 
$$
(s,t)\in S', \quad Q(s,t)=x, \quad F(s,t)=\omega. 
$$ 
Suppose further that $(s,t)$ lies in the connected component
$\mathcal{S}_1=\mathcal{S}_1(x,y,\omega,\phi)$. Now, 
if $(s_1,t_1),(s_2,t_2)\in \mathcal{S}_1$ are points such that $F(s_1,t_1)=F(s_2,t_2)$ (i.e. they lie on the same line through the origin), then we have 
$$
||(s_1,t_1)-(s_2,t_2)||\le \left(1-\frac{\sqrt{x-y}}{\sqrt{x}}\right)\cdot \max\{||(s_1,t_1)||,||(s_2,t_2)||\} \ll_{\delta}\sqrt{x}-\sqrt{x-y}\ll \frac{y}{\sqrt{x}}.
$$
If $(s_3,t_3),(s_4,t_4)\in \mathcal{S}_1$ are points such that $Q(s_3,t_3)=Q(s_4,t_4)$ (i.e., they lie on the same hyperbola $Q(\xi,\eta)=z$), then we have 
$$
||(s_3,t_3)-(s_4,t_4)||\ll_{\delta} \phi\cdot \max\{||(s_3,t_3),(s_4,t_4)||\ll_{\delta} \phi\sqrt{x}. 
$$
Hence, the distance between $(s,t)$ and any point $(\xi,\eta)\in \mathcal{S}_1$ is bounded by
$$
||(s,t)-(\xi,\eta)||\ll_{\delta} \frac{y}{\sqrt{x}}+\phi\sqrt{x}. 
$$
Using Corollary \ref{coro} and the fact that $(m,n)$ belongs to $\mathcal{M}$ if and only if the antipodal point $(-m,-n)$ belongs to $\mathcal{M}$, we deduce that there is a point $(m,n)\in \mathbb{Z}^2$ satisfying
$$
||(s,t)-(m,n)||\ll_{\delta} \frac{y}{\sqrt{x}}+\phi\sqrt{x} \ll x^{0.53/2}=Q(s,t)^{0.53/2}\ll_{\delta} ||(s,t)||^{0.53}
$$ 
such that $Q(m,n)=p$ is prime. This completes the proof of Theorem \ref{mainresult2} for the special case of $K=\mathbb{Q}(\sqrt{d})$ having class number 1 and $d\in \mathbb{N}$ being square-free with $d\equiv 2,3\bmod{4}$.  

\subsection{Proof of Theorem \ref{mainresult2} for the general case} In the case of general quadratic forms, the method is similar but the calculations become more complicated. We just indicate the required modifications below.  

The equivalence \eqref{firstcondi} remains unaltered. The function $f(\xi,\eta)$ is now defined as
$$
f(\xi,\eta)=\xi a + \eta\cdot \frac{b-\sqrt{\Delta}}{2}, 
$$
and in place of \eqref{pd2}, we have 
$$
\mathfrak{pd}=(f(m,n))
$$
if $(m,n)$ belongs to the set $\mathcal{M}$ corresponding to $\mathfrak{p}$. Now for any generator $\pi_0$ of $\mathfrak{pd}$, there is a pair $(m,n)\in \mathcal{M}$ such that $\pi_0=f(m,n)$. For an application of Theorem \ref{mainresult1}, we recall that $\mathcal{C}$ is the class whose inverse $\mathcal{C}^{-1}$ contains $\mathfrak{d}$, and hence, $\mathfrak{p}$ belongs to the class $\mathcal{C}$. Now Corollary \ref{coro} remains valid with the set $\mathcal{S}(x,y,\omega,\phi)$ defined as in \eqref{Sxyov}. In the general case, the asymptotes of the hyperbolas $Q(\xi,\eta)=x-y$ and $Q(\xi,\eta)=x$ are given by 
$$
\frac{\eta}{\xi}=\frac{2a}{-b\pm \sqrt{\Delta}}.
$$
Assume without loss of generality that $Q(\xi,\eta)$ has been transformed into reduced form and $a>0$. In this case, it can be calculated that $\mathcal{S}(x,y,\omega,\phi)$ is the intersection of the area between the two hyperbolas $Q(\xi,\eta)=x-y$ and $Q(\xi,\eta)=x$ and a sector included by the lines with slopes 
$$
s_1=-\frac{2a}{b+\sqrt{\Delta}\coth(\omega\log \epsilon)} \quad \mbox{and} \quad 
s_2=-\frac{2a}{b+\sqrt{\Delta}\coth((\omega+\phi)\log \epsilon)}.
$$
The rest of the arguments is similar as before. This completes the proof of Theorem \ref{mainresult2} in the general case.

\end{document}